\documentclass{article}[13]
\usepackage[latin1]{inputenc}

\usepackage{amssymb}
\usepackage{amsmath,enumerate,rotate}
\usepackage{amssymb,latexsym}
\usepackage{epsfig}
\usepackage{graphicx}
\usepackage[english]{babel}
\usepackage{color}

\usepackage{a4wide}

\newcommand{\RF}{{\mathcal R}}

\newcommand{\ns}{\overline{n}}

\newcommand{\PP}{{\mathbb P}}
\newcommand{\PF}{{\mathcal P}}
\newcommand{\FF}{{\mathcal F}}
\newcommand{\RR}{{\mathbb R}}
\newcommand{\NN}{{\mathbb N}}
\newcommand{\csijn}{\xi_{ij}^{(n)}}
\newcommand{\MM}{{\mathbb M}}
\newcommand{\MF}{{\mathcal M}}
\newcommand{\BF}{{\mathcal B}_1}
\renewcommand{\SS}{{\mathbb S}}
\newcommand{\SF}{{\mathcal S}}
\newcommand{\SSuno}{{\mathbb S}^{(1)}}
\newcommand{\SFuno}{{\mathcal S}^{(1)}}
\newcommand{\np}{{n'}}
\renewcommand{\ns}{{n''}}
\renewcommand{\a}{{\alpha}}

\begin{document}

\begin{center}
\bf {Central Limit Theorem with Exchangeable Summands and Mixtures of Stable Laws as Limits}
\end{center}

\vspace{1pc}
\begin{center}
Sandra Fortini\footnote{Universit\`a Bocconi, Department of Decision Sciences, via R$\ddot{o}$ntgen 1, Milano,
Italy, email: sandra.fortini@unibocconi.it}, Lucia Ladelli\footnote{Politecnico di Milano, Dipartimento di
Matematica, p.zza Leonardo da Vinci 32, 20133 Milano, Italy,  email: lucia.ladelli@polimi.it},\\ Eugenio
Regazzini\footnote{Universit\`a degli Studi di  Pavia, Dipartimento di Matematica, via Ferrata 1, 27100 Pavia,
Italy, email: eugenio.regazzini@unipv.it}
\end{center}

\begin{center}
{\em Dedicated to the memory of Enrico Magenes}
\end{center}

\vspace{3pc}






\begin{abstract}
\noindent The problem of convergence in law of normed sums of exchangeable random variables is examined. First,
the problem is studied w.r.t. arrays of exchangeable random variables, and the special role played by mixtures
of products of stable laws - as limits in law of normed sums in different rows of the array - is emphasized.
Necessary and sufficient conditions for convergence to a specific form in the above class of measures are then
given. Moreover, sufficient conditions for convergence of sums in a single row are proved. Finally, a
potentially useful variant of the formulation of the results just summarized is briefly sketched, a more
complete study of it being deferred to a future work.

\vspace{0.4cm} \noindent\textbf{AMS classification:} 60F05, 60G09

\vspace{0.4cm}

\end{abstract}

\vspace{0.4cm}

\noindent 1. {\bf Motivations and organization of the paper.} The extension of the {\em central limit theorem}
(c.l.t., for short) from independent random variables (r.v.'s, for short) to sequences of exchangeable r.v.'s
has drawn the attention of a number of researchers ever since the appearance of \cite{BlumCherRoseTeich}. Unlike
most papers, which have approached the topic in a direct way, in Section 6-7 of \cite{noi} results are derived
from statements concerning the c.l.t. for arrays of partially exchangeable r.v.'s. The consequent methods are
conducive to an original and fruitful approach to the problem. They are mentioned, for example, in a small
number of works dealing with the asymptotic behavior of the solutions of kinetic equations, because of the
recourse to a Skorokhod-type representation introduced in  \cite{noi}. See \cite{GabettaRegazziniCLT},
\cite{BassettiLadelliRegazzini}, \cite{RegazziniUMI}, \cite{GabettaRegazziniWM}, \cite{DoleraRegazzini}. With
reference to the topic developed in the present work, one can mention that the methods introduced in \cite{noi}
have been followed in \cite{RegazziniSazonov} to extend the c.l.t. to exchangeable random elements with values
in a Hilbert space. Other citations, such as in \cite{DiaconisHolmes} and \cite{Kallenberg}, are made to
complete bibliographies about exchangeability or other forms of symmetry of probability distributions (p.d.'s,
for short) of sequences of r.v.'s. We call the reader attention to \cite{JiangHahn} - a paper with which we have
been acquainted recently - because of its critical content towards our approach and consequent results. As we
are here preparing ourselves to follow such an approach to obtain new forms of the c.l.t., it is worth clearing
the field of any suspicion of mistake by recalling the recent correction note \cite{Erratum}. In it, Jiang and
Hahn, authors of \cite{JiangHahn}, admit their criticisms are mistaken and quite unjustified since ``based on a
misreading and therefore a subsequent misunderstanding of the results in \cite{noi}''.

The present paper aims at formulating conditions for the convergence of sums of exchangeable r.v.'s to random
elements distributed according to mixtures of stable laws, thus encompassing the main result in
\cite{JiangHahn}, where only the case of mixtures of Gaussians is contemplated. A short description of our
approach is contained in Section 2. Section 3 includes the formulation and the proof of the main results
concerning the weak convergence of the sequence of the sums of the elements contained in the rows of an array of
exchangeable r.v.'s. In Section 4 these results are adapted to sums of exchangeable r.v.'s. Finally, in
Section~5, a different approach - to be developed in a future paper - is mentioned. \vspace{0.4cm}

\noindent 2. {\bf Methodological background.} This section aims: $(1)$ at providing an overview of the
methodology for partially exchangeable arrays presented in \cite{noi}; $(2)$ at explaining its adaptation to the
solution of the central limit problem for sums of exchangeable r.v.'s. We begin presenting an array $A$ of {\em
exchangeable r.v.'s}, that is $A:=\{X_{ij}:\;i,j=1,2,\dots\}$. Exchangeability means that the joint distribution
of every finite subset of $m$ of these r.v.'s depends only on $m$ and not on the particular subset, $m\geq 1$.
See page 223 of \cite{ChowTeicher}. According to a classical representation theorem by de~Finetti,
exchangeability of the ${X_{ij}}'s$ is tantamount to saying that there exists a random probability measure $p^*$
such that the ${X_{ij}}'s$ turn out to be conditionally independent and identically distributed (i.i.d., for
short) given $p^*$. With symbols to be explained below, think of $p^*$ as a random element defined on $(\Omega,
{\mathcal F})$ with values in $(\PP, \PF )$. Therefore, if $\beta$ denotes the p.d. of $p^*$, assuming that all
random elements considered throughout the paper are defined on the probability space $(\Omega, {\mathcal F},
P)$, one obtains that
\begin{equation}
    \label{deFinettirappr}
P(\cap_{i=1}^{k}\cap_{j=1}^{m}\{X_{ij}\in A_{ij}\})=\int_{\PP}(\prod_{i=1}^{k}
\prod_{j=1}^{m}p(A_{ij}))\beta(dp)
\end{equation}
holds for every $A_{ij}$ in $\RF$ $(1\leq i \leq k$, $1\leq j \leq m)$, for any
$k,\;m$ in $\NN$. Note that $\RF$ stands for the Borel $\sigma$-field on $\RR$,
$\PP$ denotes the set of probability measures (p.m's, for short) on $(\RR , \RF )$, endowed with the topology of weak convergence of p.m's. Hence, $\beta$ is viewed as  p.m. on $(\PP, \PF )$, $\PF$ being the Borel $\sigma$-field generated by such a topology. The p.d. $\beta$ is usually referred to as {\em de~Finetti's measure}. From now on, weak convergence of p.m.'s will be denoted by $\Rightarrow$, while the symbols
$\overset{d}{\to}$, $\overset{d}{=}$ will designate convergence in law and equality in law, respectively.

Consider sequences $(a_n)_{n\geq1}$, $(b_n)_{n\geq1}$ of real numbers
such that
\begin{align*}
 b_n>0\;\; \mbox{for every } n,\;\;b_n \rightarrow +\infty\;\;\mbox{as}
\;\;n \rightarrow +\infty
\end{align*}
and set
\begin{align}
    \label{csi-ij}
\csijn := \frac{X_{ij}}{b_n},\;\;\;\;\; c_n:=\frac{a_n}{b_n}
\end{align}
for every $i,j,n$, and
\begin{align}
    \label{sn-cn}
S_{in}-c_n:= \sum_{j=1}^n \csijn - c_n=\frac{\sum_{j=1}^nX_{ij}-a_n}{b_n}
\end{align}
for any $i$ in $\NN$. The literature, mentioned in Section 1, deals directly with the convergence of $(S_{in}-c_n)_{n\geq1}$ for a {\em single} $i$ arbitrarily fixed (thanks to the exchangeability assumption). On the contrary, according to
\cite{noi}, our approach is that of tackling the same problem only after studying necessary and sufficient conditions for the convergence of the joint p.d. of the (infinite-dimensional) vector $(S_{1n}-c_n,S_{2n}-c_n,\dots )$. It should
be noted that, as clarified by the following example, conclusions drawn from
the latter study cannot be extended, without suitable adjustments, to the
solution of the former problem. In fact, the example below shows that the main difference regards the uniqueness of the representation of the limiting distribution. The example is drawn from \cite{JiangHahn}.

\vspace{0.4cm}

\noindent{\bf Example 1.} Let the ${X_{ij}}'s$ be i.i.d. Cauchy r.v.'s, and let
$a_n=0$, $b_n=n$. Then the (Cauchy) limiting characteristic function of $(S_{1n}-c_n)_{n\geq1}$
can be represented {\em both} as a mixture of Gaussians such as
\begin{align*}
e^{-|t|}=\int_0^{+\infty}e^{-\frac{1}{2}t^2\sigma^2}\sqrt{\frac{2}{\pi}}\frac{e^{-\frac{1}{2\sigma^2}}}{\sigma^2}d\sigma \;\;\;\;\;(t\in\RR)
\end{align*}
{\em and} as a mixture with mixing measure given by the point mass at the standard Cauchy p.d..
But, as far as the convergence of $((S_{1n}-c_n,S_{2n}-c_n,\dots ))_{n\geq1}$ is
concerned, one notices that
\begin{align*}
E\left(e^{i\sum_{k=1}^mt_kS_{kn}}\right)=\prod_{k=1}^m e^{-|t_k|}
\end{align*}
is valid for every $(t_1,\dots, t_m)$ in $\RR^m$ and $m$, $n$ in $\NN$. Hence, the
limiting (exchangeable) p.d. is presentable - in a {\em unique way} - as a mixture with de~Finetti's measure given by the above unit mass.
$\blacksquare$

The reference point for the strategy briefly sketched above is given by Theorem 2 in \cite{noi} concerning a {\em sequence of arrays} such as
\begin{align*}
A^{(n)}:=\left\{\csijn:\;\;i \in\NN,\;\;j\in\{1,\dots,n\}\right\}\;\;\;\;\;(n\in\NN)
\end{align*}
whose elements are {\em partially exchangeable}, in the sense that the joint
p.d. of every finite subset of $m_1$ elements from $(\xi_{i1}^{(n)})_{i\geq1}$, $m_2$ elements from $(\xi_{i2}^{(n)})_{i\geq1}$, $\dots$, $m_n$ elements from $(\xi_{in}^{(n)})_{i\geq1}$, depends only on $(m_1,\dots ,m_n)$ and not on the
particular subset, for any $m_1,\dots , m_n$ in $\NN_0:=\NN \cup \{0\}$. There is an extension
of de~Finetti's representation theorem according to which the above condition
of partial exchangeability is equivalent to the existence of a {\em unique}
p.m. $\beta_n$ on $(\PP^n, \PF^n)$ - with $\PF^n=\PF^{\otimes n}$ - such that

\begin{align*}
P ( \cap_{j=1}^{n} \cap_{i=1}^{m} \{ \xi_{ij}^{(n)} \in A_{ij} \})
= \int_{\PP^n}
\prod_{j=1}^n \prod_{i=1}^{m}p_j (A_{ij}) \beta_n (dp_1,\dots ,dp_n)
\end{align*}
holds for every $A_{ij}$ in $\RF$ ($1\leq j\leq n$, $1\leq i\leq m$), $n$  and $m$ in $\NN$.  In addition to this representation, partial
exchangeability can be characterized by saying that there exists a vector of random p.m.'s $(p_{1n}^*, \dots,
p_{nn}^*)$ such that the elements of $A^{(n)}$ are conditionally independent given $(p_{1n}^*, \dots, p_{nn}^*)$
and the elements of each sequence $(\csijn)_{i\geq 1}$ are conditionally identically distributed with common
p.d. $p_{jn}^*$ given $(p_{1n}^*, \dots, p_{nn}^*)$ for $j=1,\dots n$. The most general proposition in
\cite{noi} states necessary and sufficient conditions in order that, for some suitable sequence $(c_n)_{n\geq
1}$ of real numbers, the infinite-dimensional vectors $(S_{1n}-c_n,S_{2n}-c_n,\dots )$ converge in law as $n$
goes to infinity, with $S_{in}:= \sum_{j=1}^n\csijn$. Moreover, it provides a complete characterization of the
limiting p.d.'s under an assumption of {\em uniform asymptotic negligibility} for the elements of $A^{(n)}$,
which reads \vspace{0.2cm}

(H) For every $\epsilon >0$, $M_n(\epsilon)$ converges in probability to
zero as $n$ goes to infinity, where $M_n(\epsilon)$ is defined, for every $n$,
by
\begin{align*}
M_n(\epsilon)= \max_{1\leq j\leq n} p_{jn}^*([-\epsilon , \epsilon]^c).
\end{align*}
To formulate the c.l.t of interest, additional notation is needed. Let
$\overline{\RR}$ and $\overline{\RF}$ denote the extended real line and the Borel $\sigma$-field on it, respectively. In addition let $\overline{\MM}$ stand for the class of all
finite (positive) measures on $\overline{\RF}$, endowed with the topology of weak convergence.
$\overline{\MF}$ will indicate the $\sigma$-field generated by such a topology.
The symbols ${\MM}$ and ${\MF}$ will be employed with the above very same
meaning, but with regard to finite measures on $(\RR , \RF)$.
The subset of measures $\lambda$ without mass neither at $-\infty$ nor at
$+\infty$ will be denoted by $\overline{\MM}_0$, i.e.  $\overline{\MM}_0:=\{ \lambda
\in \overline{\MM}: \lambda(\{-\infty, +\infty\})=0\}$.
Now, for given $n$ in $\NN$, $j$ in $\{1,\dots , n\}$ and $\tau>0$, set
\begin{align}\label{mediatronc}
m_{jn}^{(\tau)}:= \int_{|x|\leq \tau} xp_{jn}^* (dx)
\end{align}
and consider the random p.m. $p^{*(\tau)}_{jn}$ on $\RF$ defined by
\begin{align}\label{pjn*}
p^{*(\tau)}_{jn}(A):=p^{*}_{jn}(\{x\in \RR :x- m_{jn}^{(\tau)}
\in A \})   \qquad (A\in \RF ).
\end{align}
Next, put
\begin{align}\label{mediatrasl}
\mu_n^{(\tau)}:= \sum_{j=1}^n\{ m_{jn}^{(\tau)} + \int_{\RR} \frac{x}{1+x^2}
p^{*(\tau)}_{jn}(dx) \} - c_n
\end{align}
and
\begin{align}\label{psi}
\psi_n^{(\tau)}(B) = \sum_{j=1}^n\int_{B}\frac{y^2}{1+y^2}p^{*(\tau)}_{jn}(dy)
\qquad (B\in \overline{\RF}).
\end{align}
Note that $\psi_n^{(\tau)}$ turns out to be a random measure on $\Omega$ with
values in $\overline{\MM}_0$, whereas $\mu_n^{(\tau)}$ is a real-valued r.v..
This allows us to define the p.d. $\nu_n^{(\tau)}$ of the random vector
$(\mu_n^{(\tau)} , \psi_n^{(\tau)})$ as a p.m. on $\RF \otimes
\overline{\MF}_0$, $\overline{\MF}_0$ being of course the Borel $\sigma$-field
on $\overline{\MM}_0$. Now the way is paved for the formulation of
\vspace{0.4cm}

\noindent{\bf Theorem 1.} (Theorem 2 in \cite{noi}). {\em Assume that} $A^{(n)}$
{\em is a sequence of partially exchangeable arrays with de~Finetti's
measures} $\beta_1 ,\beta_2 , \dots $, {\em meeting} (H).

\noindent{$1^o.$} {\em If} $((S_{1n} - c_n , S_{2n} - c_n , \dots ))_{n \geq 1}$
{\em converges in law, then there is a unique p.m.} $\nu$ {\em on} $(
\RR \times \overline{\MM}_0 , \RF \otimes \overline{\MF}_0) $ {\em such that
the limiting law of} $((S_{1n} - c_n , S_{2n} - c_n , \dots ))_{n \geq 1}$
{\em has Fourier-Stieltjes transform}
\begin{align}\label{theorem1}
\phi_m(t_1, \dots , t_m) = \int_{\RR \times \overline{\MM}_0} \prod_{k=1}^m
e^{\psi(t_k; \mu , \rho )} \nu (d\mu  d\rho)
\qquad ((t_1, \dots , t_m) \in \RF^m)
\end{align}
{\em for every $m$ in $\NN$, where} $e^{\psi}$ {\em is the infinitely divisible characteristic function with}
\begin{align*}
\psi(t; \mu , \rho )=i\mu t +\int_{\RR} \left(e^{itx}-1-\frac{itx}{1+x^2}
\right)\frac{1+x^2}{x^2} \rho(dx) \qquad (t\in \RF).
\end{align*}

\noindent{$2^o.$} {\em In order that} $((S_{1n} - c_n , S_{2n} - c_n , \dots ))_{n \geq 1}$
{\em converges in law - for some numerical sequence $(c_n)_{n\geq 1}$ - to a
sequence with law characterized as in \eqref{theorem1}, it is necessary and sufficient that} $\nu_n^{(\tau)}\Rightarrow \nu$ {\em for some
strictly positive} $\tau$ {\em arbitrarily fixed}.
\vspace{0.4cm}

The former part of the theorem provides a complete characterization of the family of all possible limiting laws
for the sequence $((S_{1n}-c_n,S_{2n}-c_n, \dots ))_{n\geq 1}$ of infinite-dimensional vectors of sums of
partially exchangeable r.v.'s satisfying (H). It turns out that such a family is constituted of all exchangeable
laws that are mixtures of infinitely divisible distributions. Hence, the de~Finetti measure is a p.d. for the
characteristic parameters $(\mu, \rho)$ in the L\'evy-Khintchine representation of the characteristic functions
(c.f.'s, for short). The latter part specifies, under (H), necessary and sufficient conditions in order that the
limit is a distinguished element of the class mentioned in the former one. The proof of the theorem is based on
a Skorokhod-type representation - described in details in Section 3 of \cite{noi} - we shall use to prove some
theorems in the next section. Here it is important to emphasize the {\em uniqueness} of de~Finetti's measure
$\nu$ to represent the limiting p.d. as a mixture, a property which may fail when one confines the study to the
convergence of a single sum, e.g. $(S_{1n}-c_n)_{n\geq 1}$. It is worth noting that, under such an assumption,
if the $\csijn$'s are defined as in \eqref{csi-ij} for an exchangeable array, then (H) is satisfied and the
sequence of the laws of the infinite-dimensional vector $(S_{1n}-c_n,S_{2n}-c_n,\dots )$ turns out to be
(uniformly) tight. Therefore, there exists of a subsequence $(n')$ along which the laws converge weakly to a
p.m. as in \eqref{theorem1} with a de~Finetti measure $\nu=\nu'$ which depends, in general on $(n')$. This
remark, which plays an important role in the rest of the paper, will be made precise in the next section.

The section concludes with an explanation of the main assumption made in Theorem 1 - that is partial
exchangeability - from the point of view of statistical modeling. Consider a potentially infinite sequence of
observations affected by random errors. Assume that the error for the $i$-th observation is presentable as sum
of $n$ random quantities like
\begin{align*}
S_{in}= \xi^{(n)}_{i1}+\dots + \xi^{(n)}_{in}
\end{align*}
for $i=1,2,\dots$. Suppose, in addition, that the observations are made under {\em homogeneous} conditions.
Classically, when these circumstances are taken for granted, the r.v's $\xi^{(n)}_{i1},\dots , \xi^{(n)}_{in}$
are viewed as {\em stochastically independent}, for each $i$, as well as the rows $(\xi^{(n)}_{i1},\dots ,
\xi^{(n)}_{in})$ of the array $A^{(n)}$, $i\geq 1$. In other words, the classical assumption is that the
elements of $A^{(n)}$ are stochastically independent and that the lines are identically distributed (in view of
the aforesaid homogeneity). Finally, the common p.d.'s of the $S_{in}$'s $(i=1,2,\dots)$ is approximated by the
weak limit of the law of $S_{1n}$ as $n$ goes to infinity, provided such a limit exists. The study of
conditions assuring this existence constitutes the essence of the {\em central limit problem} of probability
theory. The aforesaid independence assumption is of course unfit to deal with situations in which correlation
among different causes of error cannot be disregarded. This has led to introduce suitable forms of dependence
among the elements $\xi^{(n)}_{i1},\dots , \xi^{(n)}_{in}$, such as exchangeability. But also this assumption -
according to which the law of the causes of error would be invariant with respect to their permutations - might
be too restrictive. Therefore, making the most of the homogeneity assumption as described above, we have
preferred to consider as exchangeable only the elements of each (infinite) sequence $(\csijn)_{i\geq 1}$, for
each $j$ in $\{1, \dots , n\}$. In more precise terms, this setting can be made precise by resorting to the
scheme of partially exchangeable arrays. Then, in view of Theorem 1, the observable errors turn out to be
exchangeable according to a p.d. that can be approximated by \eqref{theorem1}, which in turn gives rise to the
Bayesian statistical model of conditionally i.i.d. observations, given $(\mu , \rho )$, with common infinitely
divisible p.d. characterized by the L\'evy-Khintchine spectral measure $\rho$ and the centering parameter $\mu$.
\vspace{0.4cm}

\noindent 3. {\bf Central limit problem for arrays of exchangeable r.v.'s.}
Going back to the same array $A:=\{X_{ij}:\;i,j=1,2,\dots\}$ of exchangeable r.v.'s as at the beginning of Section 2
and to the normed sums $S_{in}-c_n$ as in \eqref{sn-cn}, in the present section we want to deal with the following two
problems: (I) To characterize the class of all limiting laws of ${\mathbf S}_n
:= (S_{1n}-c_n , S_{2n}-c_n, \dots )$. (II) To provide a criterion (necessary
and sufficient conditions) for the weak convergence of ${\mathbf S}_n$ to a
specific form chosen from that class. As far as problem (I) is concerned, we prepare to prove that the limiting class consists of mixtures of {\em stable laws},
paralleling the well-known central limit theorem for i.i.d ${X_{ij}}'$s, in
which the admissible limiting p.d.'s are simply stable laws.
Let us recall that the c.f. of a stable law of index $\alpha \in (0,2]$ has canonical representation
\begin{equation*}
\exp\{g_{\alpha}(t; \gamma, c, \beta)\}
\end{equation*}
where
\begin{equation*}
g_{\alpha}(t; \gamma, c, \beta):= it\gamma-c|t|^{\alpha}[1+i\beta w(t, \alpha ) \text{sgn}(t)]
\end{equation*}
\begin{equation*}
  w(t,\alpha):=
  \left \{
    \begin{array}{ll}
      \text{tan} (\pi{\alpha}/2)
 &  \text{if $\alpha \neq 1$}  \\

      \text{log}|t|/{2\pi} &  \text{if $\alpha=1$ and $t\neq 0$}
    \end{array}
  \right .
\end{equation*}
with $\gamma$ in $\RR $, and $c$ in $\RR_+:=[0,+\infty)$ and $\beta$ in $B_1:=[-1,1]$. Notice that, with this parametrization, the resulting class of p.m.'s, say $\PP_s$, includes the unit masses.
Our argument is
based on Theorem 6 and the ensuing remarks in Section 6 of \cite{noi}.
For the reader's convenience, we start by providing a detailed account of that part of our previous work.
\vspace{0.4cm}

\noindent{\bf Lemma 1.} {\em In the above array $A$ assume that $(S_{1n}-c_n )_{n\geq 1}$ converges in law to a nondegenerate r.v.. Then the sequence of the laws of the ${\mathbf S}_{n}$, $n=1,2,\dots $, forms a relatively compact family of p.d.'s or, equivalently, a tight family}.
\vspace{0.4cm}

For the terminology in this statement, see Section 5 in \cite{Billingsley}. \vspace{0.4cm}

\noindent{\bf Proof.} Convergence in law of $(S_{1n}-c_n )_{n\geq 1}$, combined with exchangeability, entails
convergence of every $(S_{in}-c_n )_{n\geq 1}$, $i=1,2,\dots$. Then, designating the law of $V_n^{(\tau)}:=
(S_{1n}-c_n, \dots , S_{\tau n}-c_n)$ by $\mathcal{L}_n^{(\tau)}$, it is easy to verify that the sequence
$\mathcal{L}_n^{(\tau)}$ is (uniformly) tight. From a well-known Prokhorov's theorem, this is equivalent to
relative compactness of the same sequence: Every subsequence of $(\mathcal{L}_n^{(\tau)})_{n \geq 1}$ contains a
subsequence that converges weakly to a p.d. on $(\RR^{\tau} , \RF^{\tau} )$, and this is true for $\tau =
1,2,\dots $. Then, from Cantor's diagonal method, there exists an increasing sequence $n_1, n_2, \dots$ of
integers such that $(\mathcal{L}_{n_k}^{(\tau)})_{k\geq 1}$ converges weakly to a p.d. for every $\tau$. This is
tantamount to stating that the law of ${\mathbf S}_{n_k}$ converges weakly to a p.d. $\mathcal{L}'$, depending
on $(n'):=(n_k)_{k\geq 1}$. It is plain that the above reasoning can be repeated to state that every subsequence
of $({\mathbf S}_{n})_{n\geq 1}$ contains a subsequence which converges in law (relative compactness), which
entails uniform tightness of the laws of the ${\mathbf S}_{n}$'s in view of the Prokhorov theorem.
$\blacksquare$

We are now in a position  to formulate and prove
\vspace{0.4cm}

\noindent{\bf Theorem 2.} {\em Under the same assumptions as in} Lemma 1, {\em for every subsequence $({\mathbf S}_{n'})_{n'}$ which converges in law to an}
${\mathbf S}':=(S'_1, S'_2,\dots )$, $S'_1$ {\em being a nondegenerate r.v., there are a number} $\alpha'$ {\em in} $(0,2]$ {\em for which}
\begin{align*}
\lim_{n'\to +\infty} \frac{b_{mn'}}{b_{n'}}= m^{1/\alpha'}
\end{align*}
{\em for every} $m$ {\em in} $\NN$, {\em and a p.m} $\mu'$ {\em on the Borel} $\sigma${\em -field on}
$\Theta := \RR\times \RR_+ \times B_1$ {\em such that}
\begin{align*}
\lim_{n' \to +\infty} E\left(\mbox{exp}[i\sum_{k=1}^m t_k(S_{kn'}-c_{n'})]
\right)= \int_{\Theta} \prod_{k=1}^m \mbox{exp}\left[g_{\alpha'}(t_k; \gamma, c, \beta)\right]
\mu'(d\gamma\; dc\; d\beta )
\end{align*}
{\em holds for every} $(t_1,\dots ,t_m)$ {\em in} $\RR^m$ {\em and}
$m=1,2,\dots $, {\em with} $\mu '(d\gamma\;dc\; d\beta )\neq \delta_{\gamma_0}(\gamma)\delta_0
(dc)\mu_3(d\beta)$ {\em for every} $\gamma_0$ {\em in} $\RR$ {\em and p.m.}
$\mu_3$ {\em on} $(B_1, \BF)$, $\BF$ {\em being the Borel $\sigma$-field on $B_1$.}
\vspace{0.2cm}

\noindent{\bf Proof.} From the assumption, via Lemma~1, $V^{(m)}_{n'}:=(S_{1n'}-c_{n'}, \dots ,
S_{mn'}-c_{n'})$ converges in law to $(S'_1, \dots ,S'_m)$, for every $m$, where $(S'_k)_{k\geq 1}$
is an exchangeable sequence. Hence,
\begin{align}\label{wnm}
W_{n'}^{(m)}:= \sum_{k=1}^m(S_{kn'}-c_{n'}) \buildrel {d}\over{\to} S'_1+ \dots +S'_m
\qquad (m=1,2,\dots )
\end{align}
with $(S'_1+ \dots +S'_m)$ nondegenerate r.v. for every $m$. Then, since exchangeability entails
\begin{align}\label{tipi}
W_{n'}^{(m)} \overset{d}{=} \frac{b_{mn'}}{b_{n'}}(S_{1mn'}-c_{mn'})+\frac{a_{mn'}-ma_{n'}}{b_{n'}}
\end{align}
for every $n'$, in view of the convergence of $(S_{1n}-c_n)$ and of the convergence of types theorem (see, e.g., page 216 of \cite{Loeve}), we can say that there
are $b'_m$ and $a'_m$ such that
\begin{align*}
\frac{b_{mn'}}{b_{n'}}\to b'_m, \;\;\; \frac{a_{mn'}-ma_{n'}}{b_{n'}}\to a'_m \qquad (\np\to +\infty,\; m\in \NN).
\end{align*}
The combination of \eqref{wnm} with \eqref{tipi} gives
\begin{align}\label{limtipi}
S'_1 \buildrel {d}\over{=}  \frac{S'_1+ \dots +S'_m- a'_m}{b'_m} \qquad (m\in \NN).
\end{align}
On the other hand, exchangeability of the limiting r.v.'s $S'_i$ implies that there is a random p.m. $\pi '$
on $\RF$ with Fourier-Stieltjes transform $\tilde{\pi}'$ such that
\begin{align*}
E\left(\text{exp} \{i \sum_{k=1}^m t_k S'_k\}\right)=E\left(\prod _{k=1}^m \tilde{\pi}'(t_k)\right)
\end{align*}
holds for every $(t_1,\dots ,t_m)$ in $\RR^m$ and every $m$ in $\NN$. Whence, by resorting to
\eqref{limtipi},
\begin{align*}
E( \tilde{\pi}'(t))=E\left(\text{exp}(-it\frac{a'_m}{b'_m}) \tilde{\pi}'(\frac{t}{b'_m})^m\right) \qquad (t\in \RR ,\; m\in \NN)
\end{align*}
This and the uniqueness of de~Finetti's representation yield
\begin{align}\label{rapprstab1}
P\{\tilde{\pi}'(t)= \text{exp}(-it\frac{a'_m}{b'_m}) \tilde{\pi}'(\frac{t}{b'_m})^m \text{ for every $t$ and $m$} \}=1.
\end{align}
Hence ${\pi}'(t)$ is a stable p.d. with the exception of a set of points of $\Omega$ of $P$-probability zero.
Moreover, since the set $\PP_s$ is a closed subset of $\PP$, one gets
\begin{align}\label{rapprstab2}
E\left(\text{exp} \{i \sum_{k=1}^m t_k S'_k\}\right)=\int_{\PP_s}\left(\prod _{k=1}^m \tilde{\pi}(t_k)\right)\rho '(d\pi)
\end{align}
$\rho '$ being the restriction of the law of $\pi'$ to $\PP_s$. Now, define $T$ to be the function which
associates with each p.m. in $\PP_s$ the parameters $(\alpha, \gamma, c, \beta)$ of the canonical representation
of its c.f., with the proviso that $T(\delta_{\gamma}):=(1, \gamma,0,0)$. See, e.g., \cite{Loeve}. Then
\begin{align}\label{f.c.stabile}
\tilde{\pi}(t)= \text{exp}(g_{\alpha}(t;\gamma, c, \beta)) \qquad(t\in \RR)
\end{align}
and \eqref{rapprstab2} becomes
\begin{align*}
E(\text{exp} \{i \sum_{k=1}^m t_k S'_k\})=\int_{(0,2]\times \Theta} \prod_{k=1}^m\exp(g_{\alpha}(t_k;\gamma, c, \beta))M'(d\alpha d\gamma dc d\beta )
\end{align*}
where $M'=\rho'T^{-1}$. Now, \eqref{rapprstab1} and \eqref{f.c.stabile} give
\begin{align*}
g_{\alpha}(t;\gamma, c, \beta)= -it\frac{a'_m}{b'_m} + g_{\alpha}(\frac{t}{b'_m};\gamma, c, \beta)m \qquad (t\in \RR, \;m\in \NN)
\end{align*}
which in turn entails
\begin{equation}
\label{condizioni}
  \left \{
    \begin{array}{l}
      \gamma(1-\frac{m}{b'_m})=-\frac{a'_m}{b'_m}  \\
       c(1-\frac{m}{(b'_m)^{\alpha}})=0 \\
       c\beta[w(t, \alpha)- \frac{m}{(b'_m)^{\alpha}}w(\frac{t}{b'_m},\alpha )]=0  \qquad (m\in \NN).
    \end{array}
  \right .
\end{equation}
Recalling that $S'_1$ is a nondegenerate r.v., a consequence of \eqref{condizioni} is:
\vspace{0.2cm}

\noindent
$b'_m \equiv m \Rightarrow P(\{c=0 \}\cup \{c\neq 0, \a=\a' =1, \beta=0\})=1$ with $P(\{c=0 \}\cap \{\gamma=\gamma'\})<1$ for any $\gamma' \in \RR$;
\vspace{0.2cm}

\noindent
$b'_m \not\equiv m \Rightarrow P(\{c=0 \}\cup \{c\neq 0,b'_m \equiv m^{1/\alpha}, \alpha=\alpha'\neq 1, \gamma =-\frac{a'_m}{b'_m-m}\})=1$ with
$P(\{c=0 \})<1$.
\vspace{0.2cm}

\noindent A careful study of \eqref{condizioni} and the ensuing implications give $\a = \a'$ with $\a'$ constant and $(b'_m)_{m\geq 1}$ must coincide with $(m^{1/{\alpha'}})_{m\geq 1}$. Moreover, $M'(d\alpha\;d\gamma\;dc\;d\beta )=\delta_{\alpha'}(d\alpha)\cdot \mu'(d\gamma dc d\beta )$. $\blacksquare$
\vspace{0.2cm}

Starting from \eqref{condizioni}, further specifications about the support of $\mu'$ can be given under the same assumptions of Theorem 2.
\vspace{0.4cm}

\noindent{\bf Remark 1.} (1) If $b'_m \equiv m$, then $a'_m=0=\lim_{n'}(a_{mn'}-ma_{n'})$ for every $m$ in $\NN$ and
\vspace{0.3cm}

\noindent
$(R_1)\;\;E(\exp i\{ \sum_{k=1}^m t_k S'_k\})= \int_{\Theta}\prod_{k=1}^m \exp\{g_1(t_k;\gamma, c,0)\}\mu_{12}'(d\gamma dc)$
\vspace{0.2cm}

\noindent
$(m\in \NN,\; (t_1,\dots ,t_m) \in \RF^m)$ with either
$\mu_{12}'(\RR\times (0,+\infty))>0$ or $\mu_{12}'(d\gamma dc)=\mu_{1}'(d\gamma )\delta_0(dc)$ provided that $\mu_1'$ is a nondegenerate p.m. on $\RF$.
\vspace{0.3cm}

\noindent(2) If $b'_m\equiv  m^{1/\alpha'}$ for some $\alpha'$ in $(0,1)\cup (1,2]$, then
\vspace{0.3cm}

\noindent
$(R_2)\;\; E(\exp \{i \sum_{k=1}^m t_k S'_k\})= \int_{\RR_+\times B_1} \prod_{k=1}^m \exp\{g_{\alpha'}(t_k;\gamma', c,\beta)\}\mu_{23}'(dc d\beta)$
\vspace{0.2cm}

\noindent
$(m\in \NN,\; (t_1,\dots ,t_m) \in \RF^m)$ with $\gamma'=-a'_m/(b'_m-m)$, $\mu_{23}'((0,+\infty)\times B_1)>0$.
Note that $g_2(t;\gamma', c,\beta)=i\gamma' t -ct^2$, so that $\mu_{23}'$ can be replaced by $\mu'_2(dc)\neq \delta_0(dc)$. $\blacksquare$
\vspace{0.2cm}

A partial answer to problem (I) appears, as a straightforward consequence of Theorem~2, in the following
theorem.  The function $h$ therein, defined on $\NN$, designates, as well as in the rest of the paper, a slowly
varying function. See Appendix 1 in \cite{IbragimovLinnik}. \vspace{0.4cm}

\noindent{\bf Theorem 3.} {\em If} $({\mathbf S}_{n})=
(S_{1n}-c_n , S_{2n}-c_n, \dots )$, {\em converges in law to}
${\mathbf S}=(S_1, S_2,\dots )$, {\em then}
\begin{align}\label{cltclass}
E(\exp \{ i\sum_{k=1}^m t_k S_k\})= \int_{\Theta}\left(\prod_{k=1}^m\exp \{g_{\alpha}(t_k;\gamma,c,\beta)\}\right)\mu(d\gamma dc d\beta)
\end{align}
{\em for every $(t_1,\dots ,t_m)$ in $\RR^m$ and $m$ in $\NN$, with $\alpha$ fixed number in $(0,2]$ and $\mu$
p.m. on the Borel $\sigma$-field on $\Theta$, with the proviso that $\alpha=1$ if $\mu(\RR\times\{0\}\times
B_1)=1$. The constants $b_n$ take the form $b_n=n^{1/\alpha}h(n)$ provided that $S_1$ is nondegenerate.
Moreover, if $b_n=n h(n)$ $($i.e. $\alpha=1)$, then \eqref{cltclass} is of type $(R_1)$, while if
$b_n=n^{1/\alpha}h(n)$ with $\alpha\neq 1$ then \eqref{cltclass} is of type $(R_2)$. The mixing measure in each
representation $(R_i)$ $(i=1,2)$ is uniquely determined.} \vspace{0.4cm}

We go on to problem (II). Following \cite{noi}, we concentrate our attention on the elements involved in the
solution to the central limit problem for i.i.d. summands with common p.d. $p^*$. With a view to the convergence
criteria, these elements will be connected with problem (II) by means of a suitable application of the Skorokhod
representation theorem. In agreement with the notation introduced at the beginning of Section 2, $p^*$ denotes a
random p.m. conditionally on which the ${X_{ij}}$'s turn out to be i.i.d. with common p.d. $p^*$, and $F^*$ is
defined to be the p.d. function associated with $p^*$.  We also consider
\begin{align*}
L_n^*(x):=nF^*(xb_n)\mathbb{I}_{(-\infty,0)}(x)-n\{1-F^*(x b_n)\}\mathbb{I}_{(0,+\infty)}(x) \qquad (x \in \RR
)
\end{align*}
and the distribution function
\begin{align*}
G_n^*(x):=-nF^*(\frac{b_n}{x})\mathbb{I}_{(-\infty,0)}(x)+n\{1-F^*(\frac{b_n}{x})\}\mathbb{I}_{(0,+\infty)}(x)
\qquad (x\in\RR)
\end{align*}
along with the corresponding Lebesgue-Stieltjes measure $\lambda^*_n$ on $(\RR,\RF )$, i.e.
\begin{align*}
\lambda^*_n(A):=\int_{\RR} \mathbb{I}_{A}(x)dG_n^*(x) \qquad (A\in \RF).
\end{align*}
In addition to $L_n^*$, the so-called central convergence criterion takes into account other ``characteristic'' quantities such as
\begin{align*}
\sigma_n^*(\eta)^2:=\frac{n}{b_n^2}\left[\int_{|x|<\eta b_n} x^2p^*(dx)-\left(\int_{|x|<\eta b_n} xp^*(dx)\right)^2\right],
\end{align*}
\begin{align*}
m_n^*(\tau):=\frac{n}{b_n}\int_{|x|<\tau b_n} x p^*(dx), \;\;\;\;\;m_{1n}^*:=n\int_{\RR}\frac{b_nx}{b_n^2+x^2}p^*(dx)
\end{align*}
and
\begin{align*}
\overline{\sigma}^*(\eta)^2:=\limsup_{n\to +\infty} \sigma_n^*(\eta)^2,\end{align*}
defined for every $n$ and any positive $\tau, \; \eta$.

The elements $(m_n^*(\tau),\;m_{1n}^*)$ are seen as functions from $(\Omega, \FF)$ into $\RR^2$. Furthermore, the vector
$(\sigma_n^*(\eta)^2,\;\overline{\sigma}^*(1/n)^2)$ will be considered as a random element taking values in
$\overline\RR_+^2$, although $\sigma_n^*(\eta)^2$ is finite with probability one; this position allows to simplify future notations. As to $\lambda^*_n$, it is considered as
taking value in the space $\MM^\#$ of boundedly finite Borel (Lebesgue-Stieltjes) measures on $\RF$ endowed with
the metric
\begin{align*}
d^\#(\mu,\nu):= \int_{(0.+\infty)}\frac{d_r(\mu^{(r)},\nu^{(r)})}{1+d_r(\mu^{(r)},\nu^{(r)})}e^{-r}dr
\end{align*}
where $\mu^{(r)}$ and $\nu^{(r)}$ are the restrictions of $\mu$ and $\nu$ to $(-r,r)$, and $d_r$ is the
Prokhorov distance between the restrictions. See A2.6 in \cite{DaleyVere-Jones} for more information on
$\MM^\#$. \vspace{0.2cm}

We shall deal with the possible types of limiting distributions, singled out in Theorem~3, in the following two subsections.
\vspace{0.4cm}

\noindent 3.1 {\bf Mixtures of Gaussian distributions.} First we are interested in limiting p.d.'s for ${\bf S}_n$ characterized by
\begin{align*}
E(\exp\{i\sum_{k=1}^m  t_k S_k\})=\int_{\RR_+}\left(\prod_{k=1}^m\exp\{it_k \gamma -ct_k^2 \}\right)\pi_2(dc)
\end{align*}
for every $(t_1,\dots ,t_m)$ in $\RR^m$, $m$ in $\NN$, with $\pi_2$ different from $\delta_0$.
Gather the ``characteristic'' quantities, which are relevant to the present case, in the random vectors
\begin{align*}
W_n:=(m_n^*(\tau)-c_n, \lambda^*_n, {\sigma}_n^*(\tau)^2, p^*, \mu_n^*)
\end{align*}
\begin{align*}
W_n^{(1)}:=(p^*, \mu_n^*)
\end{align*}
for every $n$, $\mu_n^*$ being the convolution of $n$ copies of the p.d. $p^*_n$ of
$\frac{X_{ij}}{b_n}-\frac{c_n}{n}$, when $X_{ij}$ is distributed according to $p^*$ and $\tau$ is a fixed
positive number. In view of these positions, $W_n$ ($W_n^{(1)}$, respectively) turns out to be a random vector
from $(\Omega, \FF)$ into the topological product  $\SS :=\RR\times \MM^\# \times \overline\RR_+\times
\overline{\PP}_0 \times\overline{\PP}$ ($\SSuno :=\overline{\PP}_0 \times\overline{\PP}$, respectively) endowed
with the $\sigma$-field $\SF:= \RF \otimes \MF^\# \otimes \overline\RF_+\otimes \overline{\PF}_0
\otimes\overline{\PF}$ ($\SFuno: = \overline{\PF}_0 \otimes\overline{\PF}$, respectively), for every $n$. The
law of $W_n$ ($W_n^{(1)}$ respectively) will be indicated by $Q_n$ ($Q_n^{(1)}$, respectively). After denoting
the p.d. of $(m_n^*(\tau)-c_n, \lambda^*_n,{\sigma}_n^*(\tau)^2 )$ by $\nu^{(\tau)}_n$, mimicking the proof of
Lemma 1 in \cite{noi} gives \vspace{0.4cm}

\noindent{\bf Lemma 2.} {\em Let $A$ be the same array as in} Theorems 2-3.
\vspace{0.2cm}

\noindent(i) {\em If $\nu^{(\tau)}_n$ converges weakly as $n \to +\infty$ to a probability measure $\nu$ such
that $\nu(\RR\times \MM^\# \times \RR_+)=1$, then each subsequence
of $(Q_n)_{n\geq 1}$ contains a subsequence $(Q_{n'})_{n'}$, which converges weakly to a p.m. $Q'$ supported by
$\RR \times \MM^\# \times \RR_+ \times \overline{\PP}_0 \times\overline{\PP}$.}

\noindent(ii) {\em If $(S_{1n}-c_n)_{n\geq 1}$ converges in law, then each subsequence of $Q_n^{(1)}$ contains a subsequence $(Q_{n''}^{(1)})_{n''}$, which converges weakly to a p.m. $Q''^{(1)}$ supported by $\overline{\PP}_0 \times\overline{\PP}_0$.}
\vspace{0.2cm}

We now introduce a Skorokhod representation, tailored with a view to mixtures of Gaussians as limiting p.d.'s..
Since $\SS$ is a Polish space, the Skorokhod representation theorem (see, e.g., Theorem 6.7 in
\cite{Billingsley} or Theorem 11.7.2 in \cite{Dudley}) can be applied to the sequence $(Q_{n'})_{n'}$ in (i) of
Lemma 2 to state there are random elements $\hat{W}_{n'}$ and $\hat{W}$ defined on a common probability space
$(\hat{\Omega}, \hat{\FF}, \hat{P})$, such that the law of $\hat{W}_{n'}$ is $(Q_{n'})_{n'}$ for every $n'$, the
law of $\hat{W}$ is $Q'$ and $\hat{W}_{n'}(\hat{\omega})$ converges to $\hat{W}(\hat{\omega})$ w.r.t. the metric
on $\SS$, for every $\hat{\omega}\in \hat{\Omega}$.  Setting
\begin{align*}
\hat{W}_\np:=(\widehat{m_\np(\tau)-c_\np}, \hat{\lambda}_\np,  \widehat{\sigma_\np(\tau)^2}, \hat{p}, \hat{\mu}_\np)
\end{align*}
and
\begin{align*}
\hat{W}:=(\hat{m } (\tau), \hat{\lambda},  \hat{\sigma} (\tau )^2, \hat{p}, \hat{\mu})
\end{align*}
the following relations hold with the exception of a set of $\hat{P}$-probability zero:
\vspace{0.2cm}

\noindent (a) $\hat{\mu}_\np=$ convolution of $n'$ copies of $\hat{p}_\np$, where $\hat{p}_\np((-\infty,x]):= \hat{p}((-\infty, b_\np(x+\frac{c_\np}{\np})]$ for every $x$ in $\RR$.
\vspace{0.2cm}

\noindent (b) $\hat{\lambda}_\np=$ Lebesgue-Stieltjes measure associated to the distribution function
\begin{align*}
\hat{G}_\np(x):=-\np \hat{p}((-\infty,\frac{b_\np}{x}])\mathbb{I}_{(-\infty,0)}(x)+\np\hat{p}((\frac{b_\np}{x},
+\infty))\mathbb{I}_{(0, +\infty)}(x)
\end{align*}
for every $x$ in $\RR$.
\vspace{0.2cm}

\noindent (c) $\widehat{m_\np(\tau)-c_\np}=\frac{\np}{b_\np}\int_{|x|<\tau b_\np} x \hat{p}(dx)- c_\np$
\vspace{0.2cm}

\noindent (d) $\widehat{\sigma_\np(\tau)^2}=\frac{\np}{b_\np^2}\left[\int_{|x|<{b_\np\tau}} x^2 \hat{p}(dx) - \left(\int_{|x|<{b_\np \tau}} x \hat{p}(dx)\right)^2\right]$.
\vspace{0.2cm}

Passing to the sequence $(Q_{n}^{(1)})_{n}$, without any loss of generality one can state there are random elements $\hat{W}_\ns^{(1)}$ and $\hat{W}^{(1)}$ defined on the same space $(\hat{\Omega}, \hat{\FF}, \hat{P})$ such that the law of $\hat{W}_\ns$ is $Q_{n''}^{(1)}$, the law of $\hat{W}^{(1)}$ is  $Q''^{(1)}$ and $\hat{W}_\ns^{(1)}(\hat{\omega})\to \hat{W}^{(1)}(\hat{\omega})$ in the metric of $\SSuno$.
Using the same notations as in the previous point, $\hat{\mu}_\ns$ turns out to be the convolution of $\ns$ copies of $\hat{p}_\ns$ for every $\ns$, with the exception of a set of $\hat{P}$-probability zero. \vspace{0.2cm}

The above representation paves the way for a convergence criterion, that makes use of
\begin{align*}
\pi^{(\tau)}_n:&= \text{ p.d. for }\sigma^*_n(\tau)^2 \\
q^{(\epsilon)}_n:&= \text{ p.d. for }L^*_n(-\epsilon)- L_n^*(\epsilon)
\end{align*}
$\tau$ and $\epsilon$ being strictly positive numbers.
\vspace{0.4cm}

\noindent{\bf Theorem 4.} {\em In order that $({\bf S}_n)_{n\geq 1}$ derived from the exchangeable array $A$ as in} Theorem~3 {\em converge in distribution to a nondegenerate $(S_1, S_2 , \dots )$ with law characterized by}
\begin{align}\label{mixgauss}
E(\exp\{i\sum_{k=1}^m  t_k S_k\})=\int_{\RR_+}\left(\prod_{k=1}^m\exp\{it_k \gamma -\frac{1}{2}t_k^2 \sigma^2\}\right)\pi(d\sigma^2)
\end{align}
$((t_1,\dots ,t_m)$ in $\RR^m$, $m$ in $\NN)$ {\em it is necessary and sufficient that $m^*_n(\tau)-c_n
\overset{P}{\to} \gamma$, $\pi^{(\tau)}_n \Rightarrow \pi\neq \delta_0$, with $\pi(\RR_+)=1$, for some $\tau>0$
and $q^{(\epsilon)}_n\Rightarrow \delta_0$ for every $\epsilon >0$. If this is so, then $b_n=n^{1/2}h(n)$.}
\vspace{0.2cm}

\noindent{\bf Proof. Necessity.} Convergence of ${\bf S}_n$ to a nondegenerate $(S_1,S_2,\dots )$ distributed according to \eqref{mixgauss} in conjunction with Theorem 3 entail $\alpha=2$ and $b_n=n^{1/2}h(n)$. By resorting to Lemma 2 (ii) and to the ensuing Skorokhod representation, one gets
\begin{align*}
\prod_{k=1}^m \int_{\RR}e^{it_kx}\hat{\mu}_{\ns}(dx)\to \prod_{k=1}^m\int_{\RR}e^{it_kx}\hat{\mu}(dx)
\end{align*}
at each point of $\hat{\Omega}$, for every $(t_1,\dots ,t_m)$ and $m$. Hence,
\begin{align*}
\int_{\RR_+}\left(\prod_{k=1}^m\exp\{it_k \gamma -\frac{1}{2}t_k^2 \sigma^2\}\right)\pi(d\sigma^2)
=&\lim_\ns E\left(\exp \{i\sum_{k=1}^m t_k(S_{k \ns}-c_\ns)\}\right)\\
=&\lim_\ns \hat{E}(\prod_{k=1}^m \int_{\RR}e^{it_kx}\hat{\mu}_{\ns}(dx))\\
=& \hat{E}(\prod_{k=1}^m\int_{\RR}e^{it_kx}\hat{\mu}(dx))
\end{align*}
where the last equality follows from a simple application of the dominated convergence theorem. Now, in view of
the uniqueness of de Finetti's representation, there exists $\hat \sigma^2$ such that $\hat P(\hat
\sigma^2<+\infty)=1$ and
\begin{align}\label{deFinettistabile}
\int_{\RR}e^{itx}\hat{\mu}(dx)= \exp\{it\gamma-\frac{1}{2}t^2\hat{\sigma}^2 \} \qquad (t\in \RR)
\end{align}
almost surely w.r.t. $\hat{P}$. Then, from the classical central convergence criterion (cf. Subsection 23.5 of \cite{Loeve}), combined with the L\'evy representation of the Gaussian c.f. as infinitely divisible law, one gets $\hat{L}_\ns(-x)- \hat{L}_\ns(x)\to 0$ at every strictly positive $x$, which is tantamount to noting that $\hat{\lambda}_\ns$ converges to the null measure (in $\MM^\#$).
Moreover,
\begin{align*}
\hat{\sigma}_\ns(\tau )^2  \to \hat{\sigma}^2,\;\;\;
\hat{m}_\ns(\tau )- c_\ns \to \gamma.
\end{align*}
 At this stage, it should be noted that neither $\gamma$ nor the p.d. of $\hat{\sigma}^2$ depend on the specific subsequence $(\ns)$. This implies that both
$\pi^{(\tau)}_\ns$ and $q^{(\epsilon)}_\ns$ have weak limits that are independent of $(\ns)$. Therefore, by repeated application of the Skorokhod representation based on Lemma~2, the entire sequences $(m^* _n(\tau )- c_n)_{n\geq 1}$, $(\pi^{(\tau)}_n)_{n\geq 1}$
and $(q_n^{(\epsilon)})_{n\geq 1}$ converge in the sense specified by the theorem.

\noindent{\bf Sufficiency.} Let $(m^*_n(\tau )- c_n)_{n\geq 1}$, $(\pi^{(\tau)}_n)_{n\geq 1}$ and
$(q_n^{(\epsilon)})_{n\geq 1}$ be convergent in the sense specified by the theorem. Since
$q^{(\epsilon)}_n\Rightarrow \delta_0$ for every $\epsilon$, then the p.d. of $\lambda^*_n$ converges weakly to
$\delta_{\lambda_0}$, where $\lambda_0$ is the null measure. This, together with the other hypotheses, obviously
implies that the p.d. $\nu_n^{(\tau)}$ in Lemma 2 (i) converges weakly. Whence, the above Skorokhod
representation can be applied to yield
\begin{align*}
(\hat{m}_\np(\tau )- c_\np) \to \gamma,\;\;\; d^\#(\hat{\lambda}_\np,\lambda_0)\to
0,\;\;\;\hat{\sigma}^2=\lim_\np \widehat{{\sigma_\np}(\tau)^2}\end{align*} for $\np \to+\infty$ and
$\hat{\sigma}^2$ random number with p.d. $\pi$. Moreover, from the classical normal convergence criterion (see,
e.g., Section 23.5 in \cite{Loeve}), $\hat{\mu}$ belongs to $\overline{\PP}_0$ and \eqref{deFinettistabile}
holds. Finally,
\begin{align*}
\lim_\np E\left(\exp \{i\sum_{k=1}^m t_k(S_{k \np}-c_\np)\}\right)=& \lim_\np \hat{E}(\prod_{k=1}^m \int_{\RR}e^{it_kx}\hat{\mu}_{\np}(dx))\\
=&\hat{E}\left(\prod_{k=1}^m\exp\{it_k \gamma -\frac{1}{2}t_k^2 \hat{\sigma}^2\}\right)\\
=&\int_{\RR_+}\left(\prod_{k=1}^m\exp\{it_k \gamma -\frac{1}{2}t_k^2 \sigma^2\}\right)\pi(d\sigma^2)
\end{align*}
which is invariant with respect to the choice of $(\np)$ in Lemma 2 (i). Then, $((S_{1n}, S_{2n}, \dots ))$ converges in law to a sequence with law characterized by \eqref{mixgauss}. $\blacksquare$
\vspace{0.2cm}

Since a degenerate p.d. can be considered as a Gaussian law with zero variance, a simple {\em convergence criterion for degenerate limits} can be obtained from a straightforward modification of Theorem~4.
\vspace{0.4cm}

\noindent{\bf Theorem 5.} {\em In order that $({\bf S}_n)_{n\geq 1}$, derived from the exchangeable array $A$ as in} Theorems~2-3, {\em converge in distribution to $(S_1, S_2 , \dots )=(\gamma,\gamma,\dots)$,
$\gamma$ being any real number, it is necessary and sufficient that
\begin{align*}
m^*_n(\tau)-c_n  \overset{P}{\to} \gamma, \;\;\; \sigma^*_n(\tau)^2 \overset{P}{\to} 0  \;\;\;q^{(\epsilon)}_n\Rightarrow \delta_0
\end{align*}
for every $\epsilon >0$ and some $\tau>0$.}
\vspace{0.4cm}

\noindent 3.2 {\bf Mixtures of stable laws.} In agreement with Theorem~3 and Remark~1, we are now interested in the case in which the law of the limit $(S_1, S_2, \dots)$ is characterized by
\begin{align}\label{cltclass-alpha}
E(\exp \{ i\sum_{k=1}^m t_k S_k\})= \int_{\Theta}\prod_{k=1}^m\exp \{g_{\alpha}(t_k;\gamma,c,\beta)\}\mu(d\gamma dc d\beta)
\end{align}
for every $(t_1, \dots , t_m)$ in $\RR^m$, for every $m$ in $\NN$ and
for some $\alpha \neq 2$.
The ``characteristic'' quantities which pertain to the present case are gathered in the random vector
\begin{align*}
V_n:=(m_{1n}^*-c_n, \lambda^*_n, \overline{\sigma}^*(1/n)^2, p^*, \mu_n^*).
\end{align*}
It is worth taking into account, once again, the vector
\begin{align*}
W_n^{(1)}:=(p^*, \mu_n^*)
\end{align*}
$\mu_n^*$ being the convolution of $n$ copies of the p.d. $p^*_n$ of $\frac{X_{ij}}{b_n}-\frac{c_n}{n}$, when $X_{ij}$ is distributed according to $p^*$.
In view of these positions, $V_n$  turns out to be a random vector from $(\Omega, \FF)$ into the measurable space $(\SS , \SF )$ for every $n$. The law of $V_n$  will be indicated by $\tilde{Q}_n$.

After denoting the p.d. of
$(m_{1n}^*-c_n, \lambda^*_n,\bar{\sigma}^*(1/n)^2 )$ by $\kappa^{(\tau)}_n$, arguing as in the proof of Lemma~1 in \cite{noi} gives
\vspace{0.4cm}

\noindent{\bf Lemma 3.} {\em Let $A$ be the same array as in the previous theorems.}
\vspace{0.2cm}

\noindent(i) {\em If $\kappa^{(\tau)}_n$ converges weakly as $n \to +\infty$ to a probability measure $\kappa$
such that $\kappa(\RR\times \MM^\# \times \RR_+)=1$, then each
subsequence of $(\tilde{Q}_n)_{n\geq 1}$ contains a subsequence $(\tilde{Q}_{n'})_{n'}$, which converges weakly
to a p.m. $\tilde{Q}'$ supported by $\RR \times \MM^\# \times \RR_+ \times \overline{\PP}_0
\times\overline{\PP}$.}

\noindent(ii) {\em If $(S_{1n}-c_n)_{n\geq 1}$ converges in law, then each subsequence of $Q_n^{(1)}$ contains a subsequence $(Q_{n''}^{(1)})_{n''}$, which converges weakly to a p.m. $Q''^{(1)}$ supported by $\overline{\PP}_0 \times\overline{\PP}_0$.}
\vspace{0.2cm}

We describe the Skorokhod representation for the convergent subsequences $(\tilde{Q}_{n'})_{n'}$ and
$(Q_{n''}^{(1)})_{n''}$ following the same line of reasoning as in the previous subsection. Therefore, let
$\hat{V}_{n'}$ and $\hat{V}$ be random elements defined on $(\hat{\Omega}, \hat{\FF}, \hat{P})$, such that the
law of $\hat{V}_{n'}$ is $\tilde{Q}_{n'}$ for every $n'$, the law of $\hat{V}$ is $\tilde{Q}'$ and
$\hat{V}_{n'}(\hat{\omega})$ converges to $\hat{V}(\hat{\omega})$ for every $\hat{\omega}$ in $\hat\Omega$.
Setting
\begin{align*}
\hat{V}_\np:=(\widehat{m_{1\np}-c_\np}, \hat{\lambda}_\np,  \widehat{\bar{\sigma}(1/\np)^2}, \hat{p}, \hat{\mu}_\np)
\end{align*}
and
\begin{align*}
\hat{V}:=(\hat{m } , \hat{\lambda},  \hat{\sigma} ^2, \hat{p}, \hat{\mu})
\end{align*}
the following relations hold almost surely $(\hat P)$.
\vspace{0.2cm}

\noindent (a) $\hat{\mu}_\np=$ convolution of $n'$ copies of $\hat{p}_\np$ where $\hat{p}_\np((-\infty,x]):= \hat{p}((-\infty, b_\np(x+\frac{c_\np}{\np})]$ for every $x$ in $\RR$.
\vspace{0.2cm}

\noindent (b) $\hat{\lambda}_\np=$ Lebesgue-Stieltjes measure associated to the distribution function
\begin{align*}
\hat{G}_\np(x):=-\np\hat{p}((-\infty,\frac{b_\np}{x}])\mathbb{I}_{(-\infty,0)}(x)+\np \hat{p}((\frac{b_\np}{x},
+\infty))\mathbb{I}_{(0, +\infty)}(x)
\end{align*}
for every $x$ in $\RR$.
\vspace{0.2cm}

\noindent (c) $\widehat{m_{1\np}- c_\np}=\np \int_{\RR}\frac{b_\np x}{b_\np^2+x^2}\hat{p}(dx) - c_\np$
\vspace{0.2cm}

\noindent (d) $\widehat{\bar{\sigma}(1/\np)^2}=\limsup_n\frac{n}{b_n^2}\left[\int_{|x|<{\frac{b_n}{\np}}} x^2 \hat{p}(dx) - \left(\int_{|x|<{\frac{b_n}{\np}}} x \hat{p}(dx)\right)^2\right].$
\vspace{0.2cm}

Passing to the sequence $(Q_{n}^{(1)})_{n}$, one can state there are random elements $\hat{W}_\ns^{(1)}$ and $\hat{W}^{(1)}$  such that the law of $\hat{W}_\ns$ is $Q_{n''}^{(1)}$, the law of $\hat{W}^{(1)}$ is  $Q''^{(1)}$ and $\hat{W}_\ns^{(1)}(\hat{\omega})\to \hat{W}^{(1)}(\hat{\omega})$ pointwise.
Using the same notation as in Subsection~3.1, $\hat{\mu}_\ns$ turns out to be the convolution of $\ns$ copies of $\hat{p}_\ns$ for every $\ns$, with the exception of a set of $\hat{P}$-probability zero.
\vspace{0.2cm}

With a view to the next criterion, it is worth introducing further specific notation, i.e.

\vspace{0.2cm}

$\nu_{12}^{(n)}:=$ p.d. of $(m^*_{1n}-c_n, \lambda^*_n)$

$\lambda_\a:=$ Lebesgue-Stieltjes measure with distribution function $\Lambda_\a (x)= -c^-|x|^\a\mathbb{I}_{(-\infty,0)}(x)+c^+x^\a \mathbb{I}_{(0,+\infty)}(x)$, $x\in \RR$, where $c^+$ and $c^-$ satisfy $c^-, c^+ \geq 0$ and $c^-+c^+>0$
\vspace{0.2cm}

$\mathbb{G}_\a:=\{\Lambda_\a(x) :\;c^-, c^+ \geq 0,\; c^-+c^+>0\}.$
\vspace{0.2cm}

There is a one-to-one correspondence between $\mathbb{G}_\a$ and the class of all {\em L\'evy spectral measures} of stable laws of index $\alpha$. The null measure $\lambda_0$, already evoked in the proof of Theorem~4, is obtained when $c^+=c^-=0$.

\vspace{0.4cm}

\noindent{\bf Theorem 6. }{\em In order that $({\bf S}_n)_{n\geq 1}$, derived from the exchangeable array $A$ as in} Theorem~3, {\em converge in law to a nondegenerate $(S_1,S_2, \dots )$ with p.d. characterized by
\begin{align}\label{cltclass-nonuno}
E(\exp\{i\sum_{k=1}^m t_kS_k\})=\int_{\RR_+\times B_1}\left(\prod_{k=1}^m \exp\{g_\a(t_k;\gamma ,c , \beta )\}\right) \mu_{23}(dc d\beta )
\end{align}
for a fixed $\alpha$ in $(0,1)\cup (1,2)$ and for some $\mu_{23}$ with
\begin{align}\label{mixing-nonuno}
\mu_{23}((0,+\infty)\times B_1)>0
\end{align}
it is necessary and sufficient that}
\begin{align}\label{medvartronc}
&\overline{\sigma}^*(1/n)^2 \overset{P}{\to} 0,
\;\;\nu_{12}^{(n)}\Rightarrow \nu_{12} {\;\; with\;\; } \nu_{12}(\RR\times (\mathbb{G}_\a \cup\{\lambda_0\}))=1 {\;\; and\;\; } \nu_{12}(\RR\times \mathbb{G}_\a)>0.
\end{align}
{\em If this is so, then $b_n\equiv n^{1/\a}h(n)$,
\begin{align}\label{mixingalfa}
&\mu_{23}(dc d\beta)=\nu_{12}\left(\left\{(\eta, \lambda) :\frac{\pi \csc\frac{\pi\alpha}{2}}{2\Gamma(\alpha )}
\lambda(-1,1)\in dc,\; \frac{\lambda (0,1)-\lambda (-1,0)}{\lambda(-1,1)}\in d\beta \right\}\right) \\
& \!\!\!\!\!\!\! \text{with the proviso that $0/0:=0$, and} \nonumber \\
&\qquad \qquad \qquad \qquad \nu_{12}\left(\left\{(\eta, \lambda ): \eta - \frac{\lambda(0,1)-\lambda(-1,0)}{1-\a}=\gamma\right\}\right)=1
\end{align}
where $\gamma$ is the same real number as in \eqref{cltclass-nonuno}.}
\vspace{0.2cm}

\noindent{\bf Proof. Necessity.} Weak convergence of the law of ${\bf S}_n$ to a p.d. characterized by
\eqref{cltclass-nonuno}-\eqref{mixing-nonuno} for a fixed $\alpha$ in $(0,1)\cup (1,2)$ implies, through Theorem 3, that $b_n\equiv n^{1/\a}h(n)$. In view of Lemma 3 (ii) and the ensuing Skorokhod representation one gets \begin{align*}
\prod_{k=1}^m \int_{\RR}e^{it_kx}\hat{\mu}_{\ns}(dx)\to \prod_{k=1}^m\int_{\RR}e^{it_kx}\hat{\mu}(dx)
\end{align*}
and
\begin{align*}
\lim_\ns E(\exp \{ i\sum_{k=1}^m t_k (S_{k\ns}-c_\ns)\})= \int_{\RR_+\times B_1}\left(\prod_{k=1}^m\exp \{g_{\alpha}(t_k;\gamma,c,\beta)\}\right)\mu_{23}(dc d\beta).
\end{align*}
Then, from de~Finetti's representation theorem, there are $\hat{c}$ and $\hat{\beta}$ so that
\begin{align*}
\int_{\RR}e^{itx}\hat{\mu}(dx)=\exp \{g_{\alpha}(t;\gamma, \hat{c},\hat{\beta})\} \qquad (t\in \RR)
\end{align*}
with $\hat{P}$-probability one. This representation is independent of the choice of $(\ns)$ and then, by
repeated application of the Skorokhod representation based on Lemma~3, what is valid for $(\ns)$ turns out to be
valid for the entire sequence. This implies that either $\hat{\mu}$ is degenerate or $\hat{p}$ belongs to the
domain of attraction of the stable law of index $\alpha$ and parameters $\gamma,\; \hat{c},\; \hat{\beta}$.
Since, in view of \eqref{mixing-nonuno}, the latter case holds with positive $\hat{P}$ probability, then
$\hat{\lambda}_n \Rightarrow \hat{\lambda} \in \mathbb{G}_\a \cup \{\lambda_0\}$, $\widehat{m_{1n}}-c_n \to
\hat{\eta}$, $\widehat{\overline{\sigma}(1/n)^2}\to 0$ and \eqref{medvartronc}-(24) are valid. See the classical
central convergence criterion in Subsection 24.5 C of \cite{Loeve}.

\noindent{\bf Sufficiency.} When \eqref{medvartronc} is in force, Lemma 3 (i) and the Skorokhod representation yield $\hat{\lambda}_\np \Rightarrow \hat{\lambda}'\in \mathbb{G}_\a \cup \{\lambda_0\}$, $\hat{m}_{1\np}-c_\np \to \hat{\eta}'$, and $\widehat{\overline{\sigma}(1/\np)^2} \to 0.$
Consequently, the classical central convergence criterion applied to the sequence $(\np)$, gives $\hat{\mu} \in \overline{\PP}_0$ and $\int_{\RR}e^{itx}\hat{\mu}(dx)=\exp \{g_{\alpha}(t;\hat{\gamma}',\hat{c}',\hat{\beta}')\}$.
Furthermore, denoting the law of $(\hat{\gamma}',\hat{c}',\hat{\beta}')$ by $\mu$, one has
\begin{align*}
\mu(d\gamma dc d\beta) =
\end{align*}
\begin{align*}
\nu_{12}\left(\left\{ (\eta,\lambda): \eta-\frac{\lambda(0,1)-\lambda(-1,0)}{1-\a}\in d\gamma , \frac{\pi \csc\frac{\pi\alpha}{2}}{2\Gamma(\alpha )}
\lambda(-1,1)\in dc, \frac{\lambda (0,1)-\lambda (-1,0)}{\lambda(-1,1)}\in d\beta \right\}\right)
\end{align*}
which is independent of $(\np )$. Hence the entire sequence $({\bf S}_n)_n$ converges in distribution and the limiting law meets
\begin{align}\label{limite}
E(\exp \{ i\sum_{k=1}^m t_k S_k\})= \int_{\Theta}\left(\prod_{k=1}^m\exp \{g_{\alpha}(t_k;\gamma,c,\beta)\}\right)\mu(d\gamma dc d\beta)
\end{align}
with $\mu(\RR\times (0,+\infty) \times B_1)>0 $, thanks to (22). Then from Remark~1 one obtains that \eqref{cltclass-nonuno}-\eqref{mixing-nonuno} are valid. $\blacksquare$
\vspace{0.2cm}

Considering Theorem 3, it remains to analyze the case of limiting mixtures of stable laws of index $\alpha=1$, which appears when
$b_n=nh(n)$. In view of Remark 1(1), one sees that the admissible limiting c.f.'s are defined by
\begin{align}\label{cltclass-uno}
E(\exp \{ i\sum_{k=1}^m t_k S_k\})= \int_{\RR\times \RR_+}\left(\prod_{k=1}^m\exp \{g_1(t_k;\gamma,c,0)\right)\mu_{12}(d\gamma dc)
\end{align}
for every $(t_1, \dots ,t_m )$ in $\RR^m$, $m$ in $\NN$, $\mu_{12}$ being any p.d. satisfying
\begin{align}\label{mixinguno}
\mu_{12}\neq \delta_{\gamma_0}\cdot\delta_0 \qquad \text{for any $\gamma_0 \in \RR$.}
\end{align}
\vspace{0.4cm}

\noindent{\bf Theorem 7. }{\em In order that $({\bf S}_n)_{n\geq 1}$, derived from the exchangeable array $A$ as
in} Theorem 3, {\em converge in law to a nondegenerate $(S_1,S_2, \dots )$ with p.d. characterized by
\eqref{cltclass-uno}- \eqref{mixinguno} it is necessary and sufficient that} \vspace{0.4cm}
\begin{align}\label{medvartroncuno}
\overline{\sigma}^*(1/n)^2 \overset{P}{\to} 0,\;\;\;\nu_{12}^{(n)}\Rightarrow \nu_{12} \;\text{  with  } \nu_{12}(\RR\times (\mathbb{G}_1 \cup\{\lambda_0\}))=1 \text{   and  } \nu_{12}\neq \delta_{\gamma_0}\cdot \delta_{\lambda_0}.
\end{align}
{\em If this is so, then $b_n\equiv nh(n)$, and}
\begin{align}\label{mixmeasuno}
\mu_{12}(d\gamma dc)=
\nu_{12}\left(\left\{(\eta ,\lambda) :\eta  \in d\gamma,\; \frac{\pi}{2} \lambda(-1,1)\in dc \right\}\right)
\end{align}
\begin{align}\label{betazero}
\nu_{12}\left(\left\{(\eta, \lambda) :\; \lambda(0,1)-\lambda(-1,0)=0
\right\}\right)=1.
\end{align}
\vspace{0.2cm}

\noindent{\bf Proof. Necessity.} Weak convergence of the law of ${\bf S}_n$ to a p.d. characterized by
\eqref{cltclass-uno}-\eqref{mixinguno} implies, through Theorem 3, that $b_n\equiv n h(n)$. In view of Lemma 3 (ii) and the ensuing Skorokhod representation,
\begin{align*}
\prod_{k=1}^m \int_{\RR}e^{it_kx}\hat{\mu}_{\ns}(dx)\to \prod_{k=1}^m\int_{\RR}e^{it_kx}\hat{\mu}(dx)
\end{align*}
and
\begin{align*}
\lim_\ns E(\exp \{ i\sum_{k=1}^m t_k (S_{k\ns}-c_\ns)\})= \int_{\RR\times \RR_+}\left(\prod_{k=1}^m\exp \{g_{1}(t_k;\gamma,c,0)\}\right)\mu_{12}(d\gamma dc).
\end{align*}
Then, from de~Finetti's representation theorem, there are $\hat{\gamma}$ and $\hat{c}$ so that
\begin{align*}
\int_{\RR}e^{itx}\hat{\mu}(dx)=\exp \{g_1(t;\hat{\gamma}, \hat{c},0)\} \qquad (t\in \RR)
\end{align*}
with $\hat{P}$-probability one. This representation is independent of the choice of $(\ns)$ and, then, what is valid for $(\ns)$ turns out to be valid for the entire sequence. This implies that, for $\hat{\omega}$ fixed, either $\hat{\mu}$ is a unit mass or $\hat{p}$ belongs to the domain of attraction of the stable law of index $\alpha=1$ and parameters $\hat{\gamma},\; \hat{c}$, and $ \hat{\beta}=0$. Then $\hat{\lambda}_n \Rightarrow \hat{\lambda} \in \mathbb{G}_1 \cup \{\lambda_0\}$, $\widehat{m_{1n}}-c_n \to \hat{\eta}$ and $\widehat{\overline{\sigma}(1/n)^2}\to 0$.
Hence $\nu_{12}^{(n)}\Rightarrow \nu_{12} \text{  with  } \nu_{12}(\RR\times \mathbb{G}_1 \cup\{\lambda_0\})=1$ and \eqref{medvartroncuno}-(30) are valid. Furthermore, in view of
\eqref{mixinguno}, $\nu_{12}\neq \delta_{\gamma_0}\cdot \delta_{\lambda_0}$, for any $\gamma_0$ in $\RR$.

\noindent{\bf Sufficiency.} When \eqref{medvartroncuno} is in force, Lemma 3 (i) and the Skorokhod representation  yield $\hat{\lambda}_\np \Rightarrow \hat{\lambda}'\in \mathbb{G}_1 \cup \{\lambda_0\}$, $\hat{m}_{1\np}-c_\np \to \hat{\eta}'$, $\widehat{\overline{\sigma}(1/\np)^2 }\to 0$ and, consequently, the classical central convergence criterion - Subsection 24.5 C of \cite{Loeve} - applied to the subsequence $(\np)$ entails $\hat{\mu} \in \overline{\PP}_0$ and $\int_{\RR}e^{itx}\hat{\mu}(dx)=\exp \{g_{1}(t;\hat{\gamma}',\hat{c}',\hat{\beta}')\}$.
Furthermore, the law $\mu$ of $(\hat{\gamma}',\hat{c}',\hat{\beta}')$ meets
\begin{align*}
\mu(d\gamma dc d\beta)=
\end{align*}
\begin{align*}
& \nu_{12}\left(\left\{ (\eta,\lambda): \eta-(\lambda(0,1)-\lambda(-1,0))
\int_0^{\infty}\!\!\!\frac{\chi(y)-\sin y}{y^2}dy \in d\gamma, \; \frac{\pi}{2}\lambda(-1,1)\in dc,\;
\frac{\lambda (0,1)-\lambda (-1,0)}{\lambda(-1,1)}\in \!\!  d\beta \right\}\right)
\end{align*}
which is independent of $(\np )$. Therefore, the entire sequence $({\bf S}_n)_n$ converges in distribution to
$(S_1,S_2,\dots )$ and, by Remark~1, the limiting law is characterized by
\begin{align*}
E(\exp \{ i\sum_{k=1}^m t_k S_k\})= \int_{\RR\times \RR_+}\left(\prod_{k=1}^m\exp \{g_1(t_k;\gamma,c,0)\}\right)\mu_{12}(d\gamma dc)
\end{align*}
for every $(t_1, \dots ,t_m )$ and $m$. Hence, by the uniqueness of de~Finetti's measure, $\hat{\beta}'=0$ $\hat P$-a.s., and
\begin{align*}
\mu(d\gamma dc d\beta ) &=\nu_{12}\left(\left\{ (\eta,\lambda): \eta\in d\gamma, \; \frac{\pi}{2}\lambda(-1,1)\in dc\right\}\right)\cdot \delta_0 (d\beta) \\
&=\mu_{12}(d\gamma dc)\delta_0 (d\beta).
\end{align*}
Furthermore, since $\nu_{12}\neq \delta_{\gamma_0}\cdot \delta_{\lambda_0}$, $(S_1,S_2,\dots )$ is nondegenerate. Finally, an application of Theorem~3 yields $b_n=nh(n)$.
\vspace{0.4cm}

\noindent{\bf 4. Sums of exchangeable r.v.'s.} In this section we show how to use the results expounded in Section 3 to deal with the convergence of sums
\begin{equation}\label{unariga}
T_n-c_n = \frac{\sum_{i=1}^n X_i-a_n}{b_n} \qquad (n=1,2,\dots)
\end{equation}
where the $X_n$'s are exchangeable r.v.'s and
the $c_n$'s are the same as in \eqref{csi-ij}. Then, according to the notation introduced in Section 2, the law
of $(X_n)_{n\geq 1}$ obeys
\begin{align*}
P\{(X_1, \dots , X_n) \in A_1\times \dots \times A_n \}=E (\prod_{k=1}^n p^*(A_k))
\end{align*}
for every $A_1, \dots , A_n$ in $\RF$, $n$ in $\NN$.
To explain how some propositions in Section~3 can produce both new results and improvements on the existing ones, we proceed to embed the above problem into the same array $A:=\{X_{ij}:\;i,j=1,2,\dots\}$ of exchangeable r.v.'s as in the previous section. In fact, the statement of necessary and sufficient conditions in order that \eqref{unariga} converge in law can be viewed as equivalent to the problem of determining necessary and sufficient conditions for the convergence of the sum $S_{1n}-c_n$ concerning the elements of the first row of that array. If $(S_{1n}-c_n)_{n \geq 1}$ converges in law, from Lemma 1 every subsequence of $(\mathbf{S}_n)_n$ contains a subsequence $(\mathbf{S}_\np)_\np$ which converges in law. The limiting p.d., which in general depends on $(\np)$,
agrees with that given in Theorem 2.
In Section 6 of \cite{noi} an interesting case is studied where the aforesaid representation is invariant w.r.t $(\np)$, namely the case in which the limiting law is assumed to be Gaussian. An analogous circumstance happens of course if one assumes that $(S_{1n}-c_n)\overset{d}{\to} {\gamma}$, $\gamma$ being some real number. Then, Theorem~5 in Subsection 3.1 can be applied to provide an immediate and complete formulation of a {\em weak law of large numbers} partially studied in \cite{Stoica}.
\vspace{0.4cm}

\noindent{\bf Theorem 8.} {\em In order that $(T_n-c_n)_{n\geq 1}$ be convergent in probability to zero it is necessary and sufficient that the following conditions be satisfied for some $\tau >0$ and every $\epsilon >0$}.
\begin{align}\label{LDG}
& \frac{n}{b_n}\int_{[-\tau b_n,\tau b_n ]}xp^*(dx)-c_n\overset{P}{\to} 0 \\
& \frac{n}{b^2_n}\int_{[-\tau b_n,\tau b_n ]}x^2p^*(dx)-c^2_n\overset{P}{\to} 0 \\
& n(p^*((-\infty, -\epsilon b_n])+p^*((\epsilon b_n,+\infty))\overset{P}{\to} 0
\end{align}
\vspace{0.2cm}

In general, unlike the two instances just mentioned, the theorems in Subsections 3.1-3.2 fail to produce
necessary and sufficient conditions in order that the law of \eqref{unariga} be weakly convergent to a specific
p.d.. On the other hand, they can be used to state suitable sufficient conditions in a quite direct way. As an
illustration, we first consider the problem solved through Theorem 2.1 in \cite{JiangHahn} where the limiting
law of $(T_n-c_n)_{n\geq 1}$ is assumed to be a mixture of Gaussians. Since here the class of Gaussians is
thought of as including all point masses $\delta_a$  as $a$ varies in $\RR$, the above problem concerns
conditions under which the limiting c.f. of $T_n-c_n$ is presentable as
\begin{equation}\label{Rtre}
\int_{\RR \times\RR_+}\exp\{it\gamma-\frac{t^2}{2}\sigma^2\}\mu_{12}(d\gamma d\sigma^2)
\qquad (t\in \RR).
\end{equation}
With a view to the applicability of suitable results proved in the previous section we have to assume some
conditions which, together with convergence to \eqref{Rtre}, may produce the {\em invariance} of the limit
w.r.t. $(\np )$. Therefore, we shall assume
\begin{equation}\label{codegauss}
n(p^*((-\infty, -\epsilon b_n])+p^*((\epsilon b_n,+\infty))\overset{P}{\to}0 \qquad (\epsilon >0)
\end{equation}
which, besides being consistent with limiting forms studied in the previous sections, appears as a hypothesis also in Theorem 2.1 in \cite{JiangHahn}. As far as this theorem is concerned, we are now in a position to provide a direct proof of an improved formulation of it, as mentioned in \cite{Erratum}.
The argument is partially based on the proof of Theorem 8 in \cite{noi}. We take this opportunity to point out a few minor oversights.
On page 239 of \cite{noi}, lines 16 and 3 from the bottom, equations
\begin{align*}
\frac{S_{1\; n''_kr}- b_{n''_kr}}{a_{n''_kr}}=
\frac{S_{1\; n'_kr}- b_{n'_kr}}{a_{n'_kr}}\frac{a_{n'_kr}}{a_{n''_kr}}+ \frac{b_{n'_kr}-b_{n''_kr}}{a_{n''_kr}}+
\frac{X_{1 \; n'_kr+1}+ \dots + X_{1\; (n'_k+1)r}}{a_{n''_kr}}
\end{align*}
and
\begin{align*}
\frac{S_{1\; rn+s}-b_{rn+s}}{a_{rn+s}}=\frac{S_{1\; rn}-b_{rn}}{a_{rn}}
\frac{a_{rn}}{a_{rn+s}}+\frac{\sum_{i=rn+1}^{rn+s}X_i}{a_{rn+s}}+ \frac{b_{rn}-b_{rn+s}}{a_{rn+s}}
\end{align*}
should be
\begin{align*}
S_{1\; n''_kr}- \frac{b_{n''_kr}}{a_{n''_kr}}=
S_{1\; n'_kr}- \frac{b_{n'_kr}}{a_{n'_kr}}\frac{a_{n'_kr}}{a_{n''_kr}}+ \frac{b_{n'_kr}-b_{n''_kr}}{a_{n''_kr}}+
\frac{X_{1 \; n'_kr+1}+ \dots + X_{1\; (n'_k+1)r}}{a_{n''_kr}}
\end{align*}
and
\begin{align*}
S_{1\; rn+s}-\frac{b_{rn+s}}{a_{rn+s}}=S_{1\; rn}-\frac{b_{rn}}{a_{rn}}
\frac{a_{rn}}{a_{rn+s}}+\frac{\sum_{i=rn+1}^{rn+s}X_i}{a_{rn+s}}+ \frac{b_{rn}-b_{rn+s}}{a_{rn+s}},
\end{align*}
respectively.

\vspace{0.4cm}

\noindent{\bf Theorem 9.} {\em Suppose $(T_n-c_n)_{n\geq 1}$ converges in law to a nondegenerate r.v. $T$ distributed according to \eqref{Rtre}.
Then, if \eqref{codegauss} holds true, one of the following two cases takes place.

\noindent{Either}}
\begin{align}\label{unogauss}
& b_n\equiv \sqrt{n}h(n) \nonumber\\
& m^*_n(\tau)-c_n \overset{P}{\to} \gamma_0, \;\;\; \sigma^*_n(\tau)^2 \overset{d}{\to} \sigma^{*2} \;\;\; \text{for some}\; \tau>0  \\
& \gamma_0 \;\;\; \text{is a real number}, \;\; \sigma^{*2} \;\;\text{a (finite) r.v. with p.d.} \;\; \rho_1 \neq \delta_0, \;\;\mu_{12}=\delta_{\gamma_0}\cdot\rho_1 \nonumber
\end{align}
{\em or}
\begin{align}\label{duegauss}
& b_n\equiv {n}h(n) \nonumber \\
& m^*_n(\tau)-c_n \overset{d}{\to} \gamma^*, \;\;\; \sigma^*_n(\tau)^2 \overset{P}{\to} 0 \;\;\; \text{for some}\; \tau>0  \\
& \gamma^* \;\;\; \text{is an r.v. with p.d. $\rho_2 \neq \delta_a$ for every $a$ in $\RR$}, \;\; \mu_{12}=\rho_2 \cdot \delta_0. \nonumber
\end{align}

\noindent{\em Conversely if \eqref{unogauss} $($\eqref{duegauss}, respectively$)$ holds together with \eqref{codegauss}, then $(T_n-c_n)_{n\geq 1}$ converges in law to a nondegenerate r.v. $T$ with c.f. \eqref{Rtre} and $\mu_{12}=\delta_{\gamma_0}\cdot\rho_1$ $(\rho_2 \cdot \delta_0$, respectively$)$.}
\vspace{0.2cm}

\noindent{\bf Proof.} First, use the fact that convergence of $(T_n-c_n)_{n\geq 1}$ or, equivalently, of the  sum
$(S_{1n}-c_n)_{n\geq 1}$ associated with the first row of the array $A$, implies that each subsequence of $({\bf S}_n)_{n\geq 1}$ contains a subsequence $({\bf
S}_\np)_{\np}$ which converges in distribution to $(S_1',S_2', \dots )$. Combination of \eqref{codegauss} with
Theorem~1, applied to the above convergent subsequence, is enough to conclude that the law of $(S_1',S_2', \dots
)$ is a mixture of products of Gaussian p.d.'s in their wide sense (i.e. including unit masses). See, in
particular, Theorem 4 and Corollary 1 in \cite{noi}.
From this together with Theorem~2 and Remark~1 - that can be applied in view of the specific features of the elements of $A$ - one deduces that
either $b'_m\equiv \sqrt{m}$ and the law of $S'_1$ is of type \eqref{Rtre} with $\mu_{12}=\delta_{\gamma_0}\cdot\rho_1$, or $b'_m\equiv {m}$ and the law of $S'_1$ is of type \eqref{Rtre} with $\mu_{12}=\rho_2 \cdot \delta_0$. Moreover, mimicking the proof of Theorem~8 in \cite{noi} shows that either $b'_m\equiv \sqrt{m}$ for all convergent subsequences $({\bf S}_\np)_{\np}$, or $b'_m\equiv {m}$ for all convergent subsequences $({\bf S}_\np)_{\np}$. It should be noted that $ \gamma_0$, $\rho_1$, $\rho_2$ might depend on $(\np )$.
On the other hand, in the former case, the representations of the limiting c.f.'s of two convergent subsequences of $(T_n-c_n)$ must satisfy
\begin{equation*}
e^{it\gamma_0'}\int_{\RR_+}\exp\{-\frac{t^2}{2}\sigma^2\}\rho_1'(d\sigma^2)=e^{it\gamma_0}\int_{\RR_+}\exp\{-\frac{t^2}{2}\sigma^2\}\rho_1(d\sigma^2) \qquad (t\in\RR).
\end{equation*}
Since this equation implies that $\gamma_0=\gamma_0'$ and, in turn, $\rho_1=\rho_1'$ from identifiability of
scale mixtures of Gaussians, the limiting law of $(S_1',S_2', \dots )$ is independent of $(\np )$.
Then, each subsequence that converges weakly at all converges weakly to the law having finite dimensional c.f.'s
$\int_{\RR_+}\prod_{k=1}^m e^{it_k\gamma_0-t_k^2\sigma^2/2}\rho_1(d\sigma^2)$, which, thus, turns out to be the
weak limit of the ${\bf S}_n$'s. Hence, \eqref{unogauss} follows immediately from Theorem~4. Analogously, in the latter case, one shows that the limiting c.f. of $(T_n-c_n)$ is of type \eqref{Rtre} with
$\mu_{12}=\rho_2 \cdot \delta_0$, $\rho_2$ being independent $(\np )$,
and the previous argument, with Theorem~7 in the place of Theorem~4, leads to
\eqref{duegauss}. Finally, Theorems 4 and 7 can be invoked to verify sufficiency of \eqref{unogauss} and
\eqref{duegauss}, respectively. $\blacksquare$ \vspace{0.4cm}

As already pointed out, hypothesis \eqref{codegauss} turns out to be redundant if one assumes that the limiting law of  $T_n-c_n$ is a proper Gaussian p.d.. See Theorems~8-9 and Subsections~7.2-7.4 in \cite{noi}.

We go on to examine conditions for convergence to mixtures of stable laws different from the Gaussian. The leading role of assumption \eqref{codegauss} in the previous theorem is now played by

\begin{align}\label{codealpha}
&\nu_2^{(n)}, \text{ that is the law of $\lambda^*_n$, converges weakly to a p.d. $\nu_2$ such that} \\
&\nu_2(\mathbb G_{\a ,0}\cup \{\lambda_0\})=1 \text{ and } \nu_2(\mathbb G_{\a ,0})>0, \text{ with } \mathbb G_{\a ,0} :=\{\lambda\in \mathbb G_{\a}: \lambda(0,1)-\lambda(-1,0)=0\} \nonumber \\
&\text{and } \alpha \text{ fixed in } (0,1) \cup (1,2). \nonumber
\end{align}

It should be noticed that a weaker version of \eqref{codealpha} has been used in Theorem~6. A restriction is here introduced to guarantee that $\beta=0$, according to the meaning of $\beta$ pointed out by \eqref{mixingalfa}. Indeed, we are interested in limiting mixtures with c.f.
\begin{align}\label{alfa-ident}
e^{it\gamma}\int_{\RR_+}e^{-c|t|^\a}\rho(dc) \qquad(t\in \RR, \; \rho \neq \delta_0).
\end{align}
We restrict ourselves to this kind of limits since, as we shall see in a while, the uniqueness of the mixing
measure, when one considers the convergence in law of $T_n-c_n$, crucially depends on the identifiability of the
parametric family of laws appearing in the limiting mixture. \vspace{0.4cm}

\noindent{\bf Theorem 10.} {\em Let \eqref{codealpha} hold true. Then, $(T_n-c_n)_{n\geq 1}$ converges in law to a non degenerate r.v. distributed according to \eqref{alfa-ident} if and only if
\begin{align}\label{alphanonuno}
b_n\equiv n^{1/\a}h(n), \;\;\; m^*_{1n}-c_n \overset{P}{\to} \gamma, \;\;\; \overline{\sigma}^*(1/n)^2 \overset{P}{\to} 0
\end{align}
are satisfied w.r.t. the exchangeable law of} $(X_n)_{n\geq 1}$. If this is so, then
\begin{align*}
\rho(dc)= \nu_2\left(\left\{\lambda \in \mathbb G_{\a ,0}: \; \frac{\pi \csc \frac{\pi \a}{2}}{2\Gamma(\a)}\lambda(-1,1)\in dc \right\}\right)\neq \delta_0(dc).
\end{align*}

\vspace{0.2cm}
It should be noted that convergence to a stable p.d. with c.f. $t \to \exp\{i\gamma t
-c_0|t|^\a\}$, for  some $c_0$ in $(0,+\infty)$, occurs when, in the previous theorem, $\rho_3=\delta_{c_0}$.
\vspace{0.4cm}

\noindent{\bf Proof.} To prove necessity we embed, once again, the problem into the exchangeable array $A$,  and
consider $(S_{1n}-c_n)_{n\geq 1}$ in the place of $(T_n-c_n)_{n\geq 1}$. In view of the convergence in law of
$(S_{1n}-c_n)_{n\geq 1}$, every subsequence of $({\mathbf S}_n)_{n\geq 1}$ includes a subsequence $({\mathbf
S}_\np)_{\np}$ which converges in law. After a simple manipulation to reformulate \eqref{codealpha} in terms of
L\'evy-Khintchine spectral measures (see Section~6.2 in \cite{Galambos}), one sees that combination of
\eqref{codealpha} with Theorem~1 leads to conclude that the limiting law of $({\mathbf S}_\np)_{\np}$ is a
mixture of products of copies of the same stable law with index $\a$ and $\beta=0$ or, alternatively, of the
same point mass. On the other hand, exchangeability of the elements of $A$ can be considered, in conjunction
with the assumption $\nu_2(\mathbb G_{\a ,0})>0$, to obtain, through Remark~1 and uniqueness of de~Finetti's
representation, that the above limiting law is presentable as $(R_2)$ with $\a'=\a$, $\gamma$ nonrandom constant
and $\mu_{23}'=\mu_2'\cdot\delta_0$. Then, in view of \eqref{alfa-ident}, the relation
\begin{align*}
e^{it\gamma'}\int_{\RR_+}e^{-c|t|^\a}\mu_2'(dc) =e^{it\gamma}\int_{\RR_+}e^{-c|t|^\a}\rho(dc)
\end{align*}
must hold for every $t$ in $\RR$. This easily entails $\gamma'=\gamma$ and $\mu_2'= \rho$,  proving that the
limiting law of $({\mathbf S}_\np)_\np$ is independent of $(\np)$. Combining this with tightness of $({\mathbf
S}_n)_{n}$ yields the convergence in law of this sequence and, so, Theorem~6 can be applied - with
$\mu_{23}=\rho \cdot \delta_0$ and $\rho \neq \delta_0$ - to obtain \eqref{alphanonuno}. On the other hand,
Theorem~6 can also be invoked to prove that \eqref{alphanonuno} in conjunction with \eqref{codealpha} yield
convergence of $(T_n-c_n)_{n\geq 1}$ to a nondegenerate r.v. distributed according to \eqref{alfa-ident}.
$\blacksquare$ \vspace{0.4cm}

Finally, we deal with limiting p.d.'s having the form of mixtures of Cauchy p.d.'s. More specifically,  we shall
consider limiting c.f.'s like
\begin{align}\label{alfauno-ident}
e^{it\gamma}\int_{\RR_+}e^{-c|t|}\rho(dc) \qquad(t\in \RR, \; \rho \neq \delta_0).
\end{align}
The key hypothesis we shall now assume is
\begin{align}\label{codeuno}
&\nu^{(n)}, \text{ that is the law of } (m^*_{1n}-c_n, \lambda^*_n ), \text{ converges weakly to } \delta_{\gamma}\cdot \nu_2 \\
&\text{where } \nu_2(\mathbb{G}_1 \cup \{\lambda_0\})=1, \; \nu_2\neq\delta_{\lambda_0} .\nonumber
\end{align}
\vspace{0.4cm}

\noindent{\bf Theorem 11.} {\em Let \eqref{codeuno} be in force. Then, $(T_n-c_n)_{n\geq 1}$ converges in law to a nondegenerate r.v. with c.f. like \eqref{alfauno-ident} if and only if
\begin{align}\label{alphauno}
b_n = nh(n), \;\;\; \overline{\sigma}^*(1/n)^2 \overset{P}{\to} 0 \;\;\;
\end{align}
hold true w.r.t. the exchangeable law of $(X_n)_{n\geq 1}$. If this is so, then \begin{align*}
\nu_2(\{\lambda \in \mathbb{G}_1
:\lambda(0,1)-\lambda(-1,0)=0\})=1
\end{align*}
and
\begin{align*}
&\rho(dc)= \nu_2\left(\left\{\lambda \in \mathbb{G}_1 : \; \frac{\pi }{2}\lambda(-1,1)\in dc \right\}\right)\neq \delta_0(dc).
\end{align*}
\vspace{0.2cm}
}

We omit the proof since it follows from Theorem~7 with the same arguments used to prove Theorem~10.
\vspace{0.4cm}

\noindent 5. {\bf Final remarks.} We conclude the paper with a few comparative remarks concerning, in one hand,
the approach followed in the previous section and, on the other hand, an adaptation of it we would like to
develop in a future work. In fact, in the previous section we clung to the general lines fixed in \cite{noi} for
partially exchangeable r.v., even if exchangeability is the property of real interest for the present paper. The
main reason of this choice has to be ascribed to our wish to clarify in which sense the original setting is free
of the misunderstandings evoked in \cite{JiangHahn} erroneously. On the other hand, confining ourselves to
exchangeable random elements leads to limiting laws presentable as mixtures of stable laws w.r.t. their
canonical parametrization. In view of this fact, one wonders whether it is possible to substitute conditions,
which - like (\ref{medvartronc}), (\ref{medvartroncuno}), etc. - involve convergence of p.d.'s of random
measures, with conditions about convergence of random numbers. Indeed, a substitution of this kind could give
rise to more effective convergence criteria. As an illustration of the study we desire to develop, here we give
a description of the conditions one could consider to provide a version of Theorem 10 agreeing with the new viewpoint. It is worth stating beforehand that Theorems 2 and 3, together with Skorokhod
representation, will continue to play an important role. Start by fixing $\alpha$ in $(0,1)\cup(1,2)$ and
introduce conditions, that could be combined in various ways in the sequel:
\begin{align}\label{domainattr}
\frac{x^2q^*(x)}{\int_{-x}^x y^2p^*(dy)}\overset{a.s.}{\rightarrow} \frac{2-\a}{\a}  \qquad \text{as  }\; x\to
+\infty
\end{align}
\begin{align}\label{conv-c}
nq^*(b_n)\overset{d}{\to} c^*\not \equiv 0
\end{align}
\begin{align}\label{domainattrbis}
\frac{p^*((b_n, +\infty))-p^*((-\infty,-b_n))}{q^*(b_n)}\overset{P}{\to} 0\qquad \text{as  }\; n\to +\infty
\end{align}
with $q^*(x):= p^*((-\infty,-x])+ p^*((x, +\infty))$ defined for every $x>0$. Condition \eqref{domainattr}
imitates the necessary and sufficient condition in order that $p^*$ be attracted by some stable law of index
$\a$ in $(0,2]$. On the other hand, \eqref{domainattrbis} mimics the condition in order that $\beta$ in
$g_\a(t;\gamma, c, \beta)$ be zero. See, e.g., Section 25.2 B $(ii)$ in \cite{Loeve}. Finally, we take into
account condition
\begin{align}\label{conv-gamma}
m^*_{1n}-c_n \overset{P}{\rightarrow} \gamma.
\end{align}
The limiting c.f.'s we consider for the present example are the same as in (\ref{alfa-ident}). Then, a new convergence criterion
could read as follows: \emph{Under } (\ref{domainattr}) - (\ref{domainattrbis}), \emph{a necessary and
sufficient condition in order that $(T_n-c_n)_{n\geq 1}$ converge to a random number with c.f.}
(\ref{alfa-ident}) \emph{is that} (\ref{conv-gamma}) \emph{be in force}.

A possible advantage of this formulation could be that (\ref{codealpha}) is replaced by a condition in the real
field like (\ref{domainattr}). In any case, a comparative analysis of the two approaches could be made
efficacious by significant, illustrative examples. We intend to expand this aspect in conjunction with the theory developed in a future work.

\end{document}